\documentclass[11pt,leqno,twoside]{article}

\usepackage{amsfonts,amsmath,amsthm,amssymb}
\usepackage{enumerate}
\usepackage{graphics,graphicx,subfigure}

\usepackage{color}
\usepackage{todonotes}
\usepackage{cancel}
\usepackage{url}
\usepackage{hyperref}
\usepackage{makeidx}
\usepackage{showidx}
\usepackage{multicol}        
\usepackage{xspace}
\usepackage{stmaryrd}        
\usepackage{pifont}          
\usepackage{fancybox}        
\usepackage{bm}

 \usepackage{fancyhdr}
\theoremstyle{plain}
 \theoremstyle{definition}
 \newtheorem{lem}{Lemma}
 \newtheorem{defn}[lem]{Definition}
 \newtheorem{thm}[lem]{Theorem}
 \newtheorem{prop}[lem]{Proposition}
 \newtheorem{cor}[lem]{Corollary}
 \newtheorem{notn}[lem]{Notations}
 \newtheorem{pb}[lem]{Problem}
 \newtheorem{form}[lem]{Formulae}
 
 \newtheorem*{rk}{Remark}
 \newtheorem*{com}{Comment}
 \newtheorem*{ex}{Example}
 \theoremstyle{remark}

 \newcommand{\blem}{\begin{lem}}
 \newcommand{\elem}{\end{lem}}
 \newcommand{\bdefn}{\begin{defn}}
 \newcommand{\edefn}{\end{defn}}
 \newcommand{\bthm}{\begin{thm} }
 \newcommand{\ethm}{\end{thm}}
 \newcommand{\bprop}{\begin{prop}}
 \newcommand{\eprop}{\end{prop}}
 \newcommand{\bcor}{\begin{cor}}
 \newcommand{\ecor}{\end{cor}}
 \newcommand{\bnotn}{\begin{notn}}
 \newcommand{\enotn}{\end{notn}}
 \newcommand{\bpb}{\begin{pb}}
 \newcommand{\epb}{\end{pb}}
 \newcommand{\bform}{\begin{form}}
 \newcommand{\eform}{\end{form}}
 \newcommand{\brk}{\begin{rk}}
 \newcommand{\erk}{\end{rk}}
 \newcommand{\bcom}{\begin{com}}
 \newcommand{\ecom}{\end{com}}
 \newcommand{\bex}{\begin{ex}}
 \newcommand{\eex}{\end{ex}}
 \newcommand{\bpf}{\begin{proof}}
 \newcommand{\epf}{\end{proof}}

\newcommand{\bR}{\mathbb{R}}

\newcommand{\be}{\begin{equation}}
\newcommand{\ee}{\end{equation}}
\newcommand{\bal}{\begin{align}}
\newcommand{\eal}{\end{align}}
\newcommand{\ba}{\begin{align*}}
\newcommand{\ea}{\end{align*}}
\newcommand{\bmx}{\begin{matrix}}
\newcommand{\emx}{\end{matrix}}
\newcommand{\bbmx}{\begin{bmatrix}}
\newcommand{\ebmx}{\end{bmatrix}}
\newcommand{\bpmx}{\begin{pmatrix}}
\newcommand{\epmx}{\end{pmatrix}}
\newcommand{\bvmx}{\begin{vmatrix}}
\newcommand{\evmx}{\end{vmatrix}}

\newcommand{\wt}{\widetilde}
\newcommand{\f}{\frac}

\newcommand{\im}{\Rightarrow}

\newcommand{\inc}{\subseteq}

\newcommand{\setm}{\setminus}

\newcommand{\la}{\lambda}

\title{\vspace{-20mm}When is a subspace of $\ell_\infty^N$ isometrically isomorphic to $\ell_\infty^n$? \medskip\hrule height 1.2pt \vspace{-6mm}}
\author{Beata~Der\c{e}gowska\footnote{University of the National Education Commission, Krakow, Poland.
B.~D. is partially supported by The Excellent Mobility DNa.711/IDUB/EM/2025/01/00018 at the University of the National Education Commission.}, Simon Foucart\footnote{Texas A\&M University, College Station, USA. S.~F. is partially supported by the NSF (grant \#DMS-2505204).}, Barbara Lewandowska\footnote{Jagiellonian University, Krakow, Poland.
The research cooperation was funded by the program {\em Excellence Initiative – Research University} at the Jagiellonian University.}}
\date{\vspace{-6mm}\rule{100mm}{0.8pt}}

\newcommand\shorttitle{When is a subspace of $\ell_\infty^N$ isometrically isomorphic to $\ell_\infty^n$?}
\newcommand\authors{B.~Der\c{e}gowska, S.~Foucart, B.~Lewandowska}

\begin{document}
\maketitle

\vspace{-15mm}
\begin{abstract}
It is shown in this note that one can decide whether an $n$-dimensional subspace of $\ell_\infty^N$ is isometrically isomorphic to $\ell_\infty^n$ by testing a finite number of determinental inequalities.
As a byproduct,
an elementary proof is provided for the fact that
an $n$-dimensional subspace of $\ell_\infty^N$ with projection constant equal to one
must be isometrically isomorphic to $\ell_\infty^n$.
\end{abstract}



\vspace{-5mm}
\begin{center}
\rule{100mm}{0.8pt}
\end{center}


\paragraph{Prelude.}
The purpose of this note is to settle, in a testable manner, the question raised in~the~title.
To arrive at our answer,
an $n$-dimensional subspace $V$ of $\ell_\infty^N$ is better viewed as an $m$-codimensional subspace of $\ell_\infty^N$, $N = m+n$, written as $V = \{ x \in \bR^N: \langle f^1 , x \rangle = \cdots = \langle f^m , x \rangle = 0 \}$
for some linearly independent $f^1,\ldots,f^m \in\bR^N$.
In the simplest case $m=1$,
i.e., $V = \{f\}^\perp$,
it is known that $V \cong \ell_\infty^{N-1}$
if and only if $\|f\|_1 \le 2 \|f\|_\infty$.
This is a side-result of the determination
by Blatter and Cheney \cite{BlaChe},
way back in the 70's,
of a formula for the projection constant of hyperplanes in $\ell_\infty^N$---we will discuss projection constant soon.
For the next simpler case $m=2$,
an answer was given in~\cite{BarPap},
namely $V \cong \ell_\infty^{N-2}$ if and only if there exist linearly independent $f,g \in V^\perp$ and distinct indices $k \not= \ell$ such that $\|f\|_1 \le 2 |f_k|$ and $\|g\|_1 \le 2 |g_\ell|$.
The answer, however, does not directly provide a way to test whether $V$ is isometrically isomorphic to $\ell_\infty^{N-2}$.
The instantiation to the case $m=2$ of our forthcoming result (Theorem \ref{Thm??})  does.
Precisely, given linearly independent $f,g \in V^\perp$,
defining $\Delta^1,\ldots,\Delta^N \in \ell_\infty^N$ by $\Delta^k = f_k g - g_k f$,
one has  $V \cong \ell_\infty^{N-2}$ if and only if
$$
\mbox{there exist indices } k \not= \ell
\; \mbox{ such that }
\max\{ \|\Delta^k\|_1 , \| \Delta^\ell \|_1 \}
\le 2 \, |\Delta^k_\ell| \;\;(= 2 \, |\Delta^\ell_k|).
$$
Thus, it is only required to test $2\binom{N}{2}$ presumptive inequalities to settle our question.
It is important to note that the above condition is intrinsic to the space $V$,
in that it does not depend on the choice of linearly independent linear vectors $f$ and $g$ in  $V^\perp$:
e.g. if $f$ was replaced by $cf+dg$, $c \not= 0$, then each $\Delta^k$ would be replaced by $c \Delta^k$, which would not affect the status of the presumptive inequalities. 

\paragraph{Notation.}
The entries of a vector $x \in \bR^N$ are marked with a subscript, so that $x = [x_1,\ldots,x_N]^\top$.
Superscripts are reserved for indexing sequences of vectors.
For instance, a basis of the orthogonal complement $V^\perp$ of an $m$-codimensional space $V \inc \bR^N$ is written as $(f^1,\ldots,f^m)$.
In condensed form, we write
$$
F = \bbmx 
& \vline & & \vline & \\
\; f^1 & \vline & \cdots & \vline & f^m \; \\
& \vline & & \vline & 
\ebmx
\in \mathbb{R}^{N \times m}. 
$$ 
Persisting with this convention,
for a matrix $A \in \bR^{N \times m}$,
its entry located at the intersection of the $i$th row and the $j$th column is denoted by $a^j_i$,
its $j$th column is denoted by $a^j$,
and its $i$th row is denoted by $a_i$.
More generally, the row-submatrix of $A$ 
indexed by a set $S \inc \{1,\ldots,N\}$ is denoted by $A_S$.
As such, for $A,B \in \bR^{N \times m}$, 
Cauchy--Binet formula reads
$$
\det(A^\top B) = \sum_{|S|=m} \det(A_S) \det(B_S).
$$

\paragraph{Banach-Mazur distances and projection constants.}
The so-called Banach--Mazur distance between two finite-dimensional normed spaces $V$ and $W$ is defined\footnote{The finite-dimensionality is not essential---it simply ensures that the infimum is indeed attained.} as
$$
d(V,W) = \min \{ \|T\| \, \|T^{-1}\| : \; T \mbox{ is an isomorphism from } V \mbox{ to } W \}
\ge 1.
$$ 
Thus, a tautological answer is the question of the title can be: ``when $d(V,\ell_\infty^n) = 1$''.
Evidently,
this is not satisfying because there is no way (of which we are aware) of efficiently computing this Banach--Mazur distance. 
As for the projection constant of a subspace $V$ of $\ell_\infty^N$, 
it is defined\footnote{Strictly speaking, 
this quantity is the relative projection constant $\la(V,\ell_\infty^N)$ of $V$---we are making implicit use of the familiar fact that relative and absolute projection constants agree for subspaces of $\ell_\infty^N$, see e.g. \cite{Survey}.} as
$$
\la(V) = \min \{ \|P\|
:
\; P \mbox{ is a projection from $\ell_\infty^N$ onto $V$} \} \ge 1.
$$
It is well known that $\la(V) \le d(V,\ell_\infty^n)$ and here is a sketched argument for completeness:
consider a minimizing isomorphism $T: V \to \ell_\infty^{n}$;
by applying Hahn--Banach theorem componentwise,
extend it to $\wt{T} : \ell_\infty^N \to \ell_\infty^{n}$ while preserving its norm; 
then set $P := T^{-1} \wt{T}: \ell_\infty^N \to V$, which is a projection onto~$V$ (since $P(v) = T^{-1} T(v) = v$ for all $v \in V$) whose norm satisfies $\|P\| \le \| T^{-1}\| \, \|\wt{T} \|  = \| T^{-1}\| \, \|T \|= d(V,\ell_\infty^n)$;
and conclude with $\la(V) \le \|P\| \le d(V,\ell_\infty^n$).
As a result, $d(V,\ell_\infty^n)=1$ implies $\la(V)=1$.
Interestingly,
it is also known that $\la(V)=1$ conversely implies  $d(V,\ell_\infty^n)=1$,
although none of the many proofs of this result\footnote{The result brings up a possible quarrel between West and East claiming precedence:
it is often attributed to Nachbin~\cite{Nac},
although it seems to have been announced earlier by Akilov~\cite{Aki}, see the MathSciNet review MR0077897.} are elementary.
Our main result (Theorem \ref{Thm??}) actually provides an elementary proof of the equivalence
$\la(V) = 1 \iff d(V,\ell_\infty^n)=1$,
albeit with the restriction that $V$ is (isometrically isomorphic to) a subspace of $\ell_\infty^N$.
Thus, a better answer to our question is: ``when $\la(V) = 1$''.
Arguably,
this is a satisfying answer
because the projection constant of a subspace of $\ell_\infty^N$ can be computed by linear programming (see e.g. \cite{Survey} for details)...
except that most optimization solvers do not work in exact arithmetic,
so truly testing the equality $\la(V) = 1$ could be problematic. 
In this sense, the answer we give to the question of the title is ``more'' satisfying---it entails verifying a finite (but possibly large) number of inequalities which can,
on the face of it,
be handled symbolically.

\paragraph{The main result.}
Without further ado,
our awaited answer to the question 
``when is a subspace $V$ of $\ell_\infty^N$ isometrically isomorphic to $\ell_\infty^n$''
materializes as item (i) of the theorem below.
Its statement involves
an intrinsic basis $(h(S)^k, k \in S)$ of $V^\perp$
associated with a set $S \inc \{1,\ldots,N \}$ of size $m = {\rm codim}(V)$.
Although it is constructed by invoking a fixed basis $(f^1,\ldots,f^m)$ of $V^\perp$, 
note that it is actually independent of this fixed basis.
Its defining formula is,
for $k \in S$ and $i = 1,\ldots,N$,
$$
h(S)^k_i =  \f{\det(F_S[{\rm row}_k \leftarrow {\rm row}_i])}{\det(F_S)},
\qquad \quad \mbox{where} \quad
F = \bbmx 
& \vline & & \vline & \\
\; f^1 & \vline & \cdots & \vline & f^m \; \\
& \vline & & \vline & 
\ebmx
\in \mathbb{R}^{N \times m}. 
$$
On the one hand, the fact that the $h(S)^k$, $k \in S$, belong to $V^\perp$ follows from a Laplace expansion with respect to the $k$th row, yielding
$$
h(S)^k_i = \f{1}{\det(F_S)} \sum_{j=1}^m (-1)^{k+j} f^j_i \det(F_{S \setm \{k\}}^{[1:m] \setm \{j\}})\qquad \mbox{for all } i = 1,\ldots, N.
$$
In the particular case $m = 2$ and $S = \{k,\ell\}$,
we observe that $h(S)^k = \Delta^\ell/\Delta_k^\ell$,
which leads to the result mentioned in the prelude.
On the other hand, the fact that the $h(S)^k$, $k \in S$, are linearly independent follows from
$$
h(S)^k_i = 
\begin{cases}
0 & \mbox{ if } i \in S \mbox{ is different from } k,\\
1 & \mbox{ if } i \in S \mbox{ is identical with } k.
\end{cases}
$$
As a consequence,
any $f \in V^\perp$ is expressed as
$f = \sum_{k \in S} {f_k} h(S)^k$.
In matrix form, this can simply be written as 
the identity $F =  H(S) F_S$, to be used later.

\bthm
\label{Thm??}
Given an $m$-codimensional subspace $V$ of $\ell_\infty^N$,
the following statements are equivalent:\vspace{-5mm}
\begin{enumerate}
\item[(i)] there exists an index set $S$ of size $m$
such that $\|h(S)^k\|_1 \le 2$ for all $k \in S$;
\vspace{-1mm}
\item[(ii)] $V$ is isometrically isomorphic to $\ell_\infty^{n}$, $n = N-m$,
i.e., $d(V,\ell_\infty^n) = 1$;\vspace{-1mm}
\item[(iii)] the projection constant of $V$ equals one, i.e., $\la(V)=1$.
\end{enumerate}
\ethm

The justification of these equivalences owes to the lemmas below.
Indeed, the implication (i)$\im$(ii) follows from Lemma \ref{Lem1}, which is relatively straightforward;
the implication (ii)$\im$(iii) is a consequence of $\la(V) \le d(V,\ell_\infty^n)$;
and the implication (iii)$\im$(i) follows from Lemma \ref{Lem2},
which is the centerpiece of this note.

\blem
\label{Lem1}
For any index set $S$ of size $m$ such that $\det(F_S) \not= 0$,
$$
d(V,\ell_\infty^n) \le \max \big\{ 1, \max_{k \in S} \|h(S)^k \|_1 - 1 \big\}.
$$
\elem

\blem
\label{Lem2}
Let $P$ be (the matrix of) a projection from $\ell_\infty^N$ onto $V$ with $\|P\| = \la(V)$.
For any index set $S$ of size $m$ such that $\det(F_S) \not= 0$ and
$\det(I - P_S^S) \not= 0$,
$$
\max_{k \in S} \| h(S)^k \|_1 - 1
\le 1 + (\la(V) - 1)  
\|(I -P_S^S)^{-1}\|.
$$ 

\elem

\bpf[Proof of Lemma \ref{Lem1}]
For $v \in V = \{f^1,\ldots,f^m \}^\perp$,
the equality $F^\top v = 0$
yields
$F_S^\top v_S + F_{S^c}^\top v_{S^c} = 0$,
i.e., 
$v_S = - F_S^{- \top} F_{S^c}^\top v_{S^c}$. 
This implies that
$$
\| v_S \|_\infty  \le  \|F_S^{- \top} F_{S^c}^\top\| \, \|v_{S^c}\|_\infty,
$$ 
where the operator norm is transformed into
\begin{align*}
\|F_S^{- \top} F_{S^c}^\top\|
& = 
\max_{k \in S} \sum_{i \in S^c} \big| (F_S^{- \top} F_{S^c}^\top)_k^i \big|
= \max_{k \in S} \sum_{i \in S^c} \big| ( F_{S^c} F_S^{- 1})^k_i \big|
= \max_{k \in S} \sum_{i \in S^c} 
\Big| \sum_{j=1}^m (F_{S^c})^j_i (F_S^{- 1})^k_j  \Big|\\
& = \max_{k \in S} \sum_{i \in S^c} \bigg| \sum_{j=1}^m f_i^j \f{(-1)^{k+j} \det(F_{S \setm \{k\}}^{[1:m] \setm \{j\}})  }{\det(F_S)}  \bigg|
= \max_{k \in S} \sum_{i \in S^c} \big|  h(S)^k_i \big|
= \max_{k \in S} \|h(S)^k_{S^c} \|_1.
\end{align*}
It follows that,
for any $v \in V$,
$$
\|v\|_\infty
= \max\{ \|v_{S^c}\|_\infty, \|v_{S}\|_\infty \}
\le \max \Big\{  1 , \max_{k \in S} \| h(S)^k_{S^c} \|_1 \Big\} \|v_{S^c}\|_\infty.
$$
Since $\|v_{S^c}\|_\infty \le \|v\|_\infty$ also holds for any $v \in V$, we deduce that
$$
d(V,\ell_\infty^{N-m}) \le \max \Big\{  1 , \max_{k \in S} \|h(S)^{k}_{S^c} \|_1 \Big\}.
$$
The announced form of the result makes use $\| h(S)^{k}_{S^c} \|_1 = \|h(S)^{k} \|_1 - \|h(S)^{k}_{S} \|_1 = \|h(S)^{k} \|_1 - 1$.
\epf

\bpf[Proof of Lemma \ref{Lem2}]
Let $P$ be a (minimal) projection from $\ell_\infty^N$ onto~$V$.
Since $I-P$ vanishes on $V = \{f^1,\ldots,f^m\}^\perp$,
there exist $y^1,\ldots,y^m \in \bR^N$ 
such that $(I-P)x = \sum_{i=1}^m \langle f^i, x \rangle y^i$ for all $x \in \bR^N$.
Then, in view of $Px \in V $ for all $x \in \bR^N$,
we have $0 = \langle f^j, Px \rangle = \langle f^j,x \rangle - \sum_{i=1}^m \langle f^i, x \rangle \langle f^j,y^i \rangle$ for all $j=1,\ldots,m$.
This forces $\langle f^j, y^i \rangle = \delta_{i,j} $ for all $i,j=1,\ldots,m$.
All in all, the projection $P$ can be expressed, for any $x \in \bR^N$, as
$$
Px = x - \sum_{i=1}^m \langle f^i, x \rangle y^i,
\qquad \mbox{ where } y^1,\ldots,y^m \in \bR^N
\mbox{ satisfy }
\langle f^j, y^i \rangle = \delta_{i,j}.
$$
In a more condensed matrix form, this reads 
$$
P = I_N - Y F^\top
\qquad \mbox{ where } Y \in \bR^{N \times m} \mbox{ satisfies } F^\top Y = I_m .
$$
Relatively to another basis $(g^1,\ldots,g^m)$ of $V^\perp$,
written as $G = F M$ for some invertible matrix $M \in \bR^{m \times m}$,  we have
$$
P = I_N - Z G^\top
\qquad \mbox{ where } Z = Y M^{- \top} \in \bR^{N \times m} \mbox{ satisfies } G^\top Z = I_m .
$$
In view of $\sum_{|S|=m} \det(F_S) \det(Y_S) = 1$,
which stems from Cauchy--Binet formula,
we can find an index set $S$ such that not only $\det(F_S) \not= 0$ but also $\det(Y_S) \not= 0$.
The former is needed in the definition of the $h(S)^k$, $k \in S$, and the latter will be needed soon.
Fixing this index set $S$ from now~on,
we take $(g^1,\ldots,g^m)$ to be the basis $(h^k, k \in S)$---dropping the dependence on $S$ for ease of notation.
The matrices $G$, $Z$, and $P$ thus take the form
$$
H = \bbmx I_m \\ \hline H_{S^c} \ebmx,
\qquad 
Z = \bbmx Z_S \\ \hline Z_{S^c} \ebmx,
\qquad 
P 
= I_N - \bbmx
Z_S & \vline & Z_S H_{S^c}^\top \\
\hline
Z_{S^c} & \vline & Z_{S^c} H_{S^c}^\top 
\ebmx.
$$
From this expression of $P$,
it follows that 
\begin{align*}
\|P\| & = \max_{i=1,\ldots,N} \sum_{j=1}^N |P^j_i|
\ge \max_{i \in S} \bigg(
|1-Z_{i}^i| + \sum_{j \in S \setm \{i\}} |Z_i^j| 
+ \sum_{j \in S^c} |(Z_S H_{S^c}^\top)_i^j|
\bigg)\\
& \ge \max_{i \in S} \bigg(
1- |Z_{i}^i| + \sum_{j \in S \setm \{i\}} |Z_i^j| 
+ \sum_{j \in S^c} |(Z_S H_{S^c}^\top)_i^j|
\bigg).
\end{align*}
Therefore, for any $i \in S$,
we obtain after some rearrangement
$$
\|P\| - 1 + \alpha_i \ge \beta_i,
\qquad \mbox{where }
\alpha_i:=  |Z_{i}^i| - \sum_{j \in S \setm \{i\}} |Z_i^j|
\quad \mbox{and} \quad
\beta_i := \sum_{j \in S^c} |(Z_S H_{S^c}^\top)_i^j|.
$$
For any $c \in \bR^S$,
we observe on the one hand that 
$$
\sum_{i \in S} \beta_i |c_i| 
= \sum_{j \in S^c} \sum_{i \in S} |(Z_S H_{S^c}^\top)_i^j| \, |c_i|
\ge \sum_{j \in S^c}\Big| \sum_{i \in S} 
(H_{S^c} Z_S^\top)_j^i c_i \Big|
= \sum_{j \in S^c} \big| (H_{S^c} Z_S^\top c)_j \big|,
$$
and on the other hand that
\begin{align*}
\sum_{i \in S} \alpha_i |c_i|
& = \sum_{i \in S} |Z_{i}^i| \, |c_i| - \sum_{\substack{i,j \in S \\ i \not= j}} |Z_i^j| \, |c_i|
= \sum_{j \in S} |Z_{j}^j| \, |c_j| - \sum_{\substack{i,j \in S \\ i \not= j}} |Z_i^j| \, |c_i|\\
& = \sum_{j \in S} \Big(
|Z_{j}^j| \, |c_j| - \sum_{i \in S \setm \{ j \} } |Z_i^j| \, |c_i|
\Big)
\le \sum_{j \in S} \Big| \sum_{i \in S} Z_i^j c_i \Big|
= \sum_{j \in S} \big| (Z_S^\top c)_j \big|.
\end{align*}
We consequently derive that, for any $c \in \bR^S$,
$$
\big( \|P\| - 1 \big) \sum_{i \in S} |c_i|
+ \sum_{j \in S} \big| (Z_S^\top c)_j \big| \ge 
\sum_{j \in S^c} \big| (H_{S^c} Z_S^\top c)_j \big|.
$$
At this point,
we need the specificity of the index set $S$ to ensure that the matrix $Z_S$ is invertible.
This holds true thanks to the identity $F = H F_S$,
i.e.,  $H = F M$ with $M = F_S^{-1}$,
which implies that $Z = Y M^{-\top} = Y F_S^\top$,
so $Z_S = Y_S F_S^\top$ is invertible as the product of two invertible matrices.
Hence, for any $\ell \in S$,
we can make the choice $c = Z_S^{-\top} h_S^\ell$,  
for which $c_i = (Z_S^{-1})_\ell^i$ and
 $Z_S^\top c = h_S^\ell = \delta^\ell$,
to arrive at
$$
\big( \|P\| - 1 \big) \sum_{i \in S} \big| (Z_S^{-1})_\ell^i \big|
+ 1 \ge 
\sum_{j \in S^c} \big| h^\ell_j \big|.
$$
Restoring the dependence on $S$, we have shown that there exists an index set $S$
(any $S$ such that $\det(F_S)\not= 0$ and $\det(Y_S) \not= 0$ is suitable)
such that
$$
\max_{\ell \in S} \|h(S)^\ell_{S^c}\|_1
\le 1 + \big( \|P\| - 1 \big) 
\max_{\ell \in S} \sum_{i \in S} \big| (Z_S^{-1})_\ell^i \big|.
$$ 
Taking into account that $\|P\| = \la(V)$ for a minimal projection,
recognizing that the last maximum is $\|Z_S^{-1}\|$,
and identifying $Z_S$ with $I - P_S^S$,
as apparent from the block-representation of $P$, completes the argument.
\epf


\begin{thebibliography}{99}

\bibitem{Aki}
Akilov, G. P. (1947).
{\em On the extension of linear operations.} 
Doklady Akad. Nauk SSSR (N.S.), 57, 643--646.

\bibitem{BarPap}
Baronti, M., Papini, P. (1991). 
{\em Norm-one projections onto subspaces of finite codimension in $\ell_1$ and $c_0$.} 
Periodica Mathematica Hungarica, 22, 161--174.

\bibitem{BlaChe} 
Blatter, J., Cheney, E. W. (1974). 
{\em Minimal projections on hyperplanes in sequence spaces.} Annali di Matematica Pura ed Applicata, 101, 215--227.

\bibitem{Survey}
Foucart, S., Skrzypek, L. (202x).
{\em Minimal projections: from classical theory to modern developments.}
Surveys in Approximation Theory.
In preparation.

\bibitem{Nac}
Nachbin, L. (1950). 
{\em A theorem of the Hahn--Banach type for linear transformations.} 
Transactions of the American Mathematical Society, 68(1), 28--46.

\end{thebibliography}
\end{document}